\input amstex
 
\input amsppt.sty
\hsize = 6.5 truein
\vsize = 9 truein
\NoBlackBoxes
\TagsAsMath
 
\newskip\sectionskipamount
\sectionskipamount = 24pt plus 8pt minus 8pt
\def\sectionskip{\vskip\sectionskipamount}
\define\sectionbreak{%
        \par  \ifdim\lastskip<\sectionskipamount
        \removelastskip  \penalty-2000  \sectionskip  \fi}
\define\section#1{%
        \sectionbreak                   
        \subheading{#1}%
        \bigskip
        }
 
 
        
        \let    \< = \langle
        \let    \> = \rangle
 
\define\op#1{\operatorname{\fam=0\tenrm{#1}}} 
\define\opwl#1{\operatornamewithlimits{\fam=0\tenrm{#1}}}
 
        \define         \a              {\alpha}
        \redefine       \b              {\beta}
        \redefine       \d              {\delta}
        \redefine       \D              {\Delta}
        \define         \e              {\varepsilon}
        \define         \g              {\gamma}
        \define         \G              {\Gamma}
        \redefine       \l              {\lambda}
        \redefine       \L              {\Lambda}
        \define         \n              {\nabla}
        \redefine       \var    {\varphi}
        \define         \s              {\sigma}
        \redefine       \Sig    {\Sigma}
        \redefine       \t              {\tau}
        \define         \th             {\theta}
        \redefine       \O              {\Omega}
        \redefine       \o              {\omega}
        \define         \z              {\zeta}
        
        \redefine       \i              {\infty}
        \define         \p              {\partial}
 
\def\prd{\Pi}

\topmatter
\title On a Conjecture of Tarski on Products of Cardinals
\endtitle
\author Thomas Jech${}^1$ and Saharon Shelah${}^2$
\endauthor
\address{%
T\. Jech, Department of Mathematics, The Pennsylvania State University, 
University Park, PA~~16802 \newline
S\. Shelah, Department of Mathematics, The Hebrew University, Jerusalem, 
Israel}\endaddress
\thanks{${}^1$Supported partially by an NSF grant.  I wish to express my 
gratitude to the Mathematical Institute of the Eidgen\"ossische Technische 
Hochschule in Z\"urich for its hospitality during my visit.\newline
${}^2$Publication \#385.  Supported partially by the B\.S\.F.\newline
${}^3$AMS classification:  03E\newline
${}^*$Keywords:  Cardinal arithmetic, singular cardinals problem, $pcf$}
\endthanks
\endtopmatter

\document
 
\sectionskip
 
\heading Abstract${}^{3*}$
\endheading
 
\vskip 0.5 truein
 
We look at an old conjecture of A\.~Tarski on cardinal arithmetic and show 
that if a counterexample exists, then there exists one of length $\o_1 + \o$.
 
\newpage
In the early days of set theory, Hausdorff and Tarski established basic rules 
for exponentiation of cardinal numbers.  In \cite{T} Tarski showed that for 
every limit ordinal $\b$, $\prod_{\xi < \b} \aleph_{\xi} = 
\aleph_{\b}{}^{|\b|}$, and conjectured that
$$
\prod_{\xi < \b} \aleph_{\s_{\xi}} = \aleph_{\a}{}^{|\b|} \tag 1
$$
holds for every ordinal $\b$ and every increasing sequence $\{\s_{\xi}\}_{\xi 
< \b}$ such that $\lim_{\xi < \b} \s_{\xi} = \a$.  He remarked that $(1)$ 
holds for every countable ordinal $\b$.
 
\demo{Remarks} 1. The left hand side of (1) is less than or equal to
the right hand side.
 
2. If $\b$ has $|\b|$ disjoint cofinal subsets then the equality (1) holds.
Thus the first limit ordinal that can be the length of a counterexample
to (1) is $\o_1 + \o.$

\enddemo
 
\demo{[Proof} Let $\{A_i : i < |\b|\}$ be disjoint cofinal subsets of 
$\b.$ Then \newline $\prod_{\xi < \b} \aleph_{\s_{\xi}} \ge
\prod_{i < |\b|} \prod_{\xi \in A_i} \aleph_{\s_{\xi}} \ge 
\prod_{i < |\b|} \aleph_{\a} = \aleph_{\a}{}^{|\b|}.]$
\enddemo
 
It is not difficult to see that if one assumes the Singular Cardinals 
Hypothesis then $(1)$ holds.  With the hindsight given by results obtained in 
the last twenty years, it is also not difficult to find a counterexample to 
Tarski's conjecture.  For instance, using the model described in \cite{M}, one 
can have an increasing sequence of cardinals of length $\b = \o_1 + \o$ whose 
product does not satisfy $(1)$.  The purpose of this note is to show that 
if Tarski's conjecture fails then it fails in this specific way. Namely,
if there is a counterexample then there is one of length $\o_1 + \o.$
 
The main result of this paper is the following:
 
\proclaim{Theorem} A necessary and sufficient condition for Tarski's 
conjecture to fail is the existence of a singular cardinal $\aleph_{\g}$ of 
cofinality $\aleph_1$ such that $\aleph_{\g} > \aleph_{\o_1}{}^{\aleph_1}$ and 
$\aleph_{\g}{}^{\aleph_1} > \aleph_{\g+\o}{}^{\aleph_0}$.
\endproclaim
 
If $\aleph_{\g}$ is a cardinal that satisfies the condition then the sequence 
$\{\aleph_{\xi}\}_{\xi < \o_1} \cup \{\aleph_{\g+n}\}_{n < \o}$ is a 
counterexample to $(1)$:
$$
\prod_{\xi < \o_1} \aleph_{\xi} \cdot \prod_{n < \o} \aleph_{\g+n} = 
\aleph_{\o_1}{}^{\aleph_1} \cdot \aleph_{\g+\o}{}^{\aleph_0} < 
\aleph_{\g+\o}{}^{|\o_1+\o|}.
$$
Such a cardinal exists in one of Magidor's models, e\.g\. when $\aleph_{\g} = 
\aleph_{\o_1+\o_1}$ is a strong limit, $\aleph_{\o_1+\o_1}{}^{\aleph_1} = 
\aleph_{\o_1+\o_1+\o+2}$ and $\aleph_{\o_1+\o_1+\o}{}^{\aleph_0} = 
\aleph_{\o_1+\o_1+\o+1}$.  
 
Also, if $\l > \aleph_{\o_1}$ is a strong limit 
singular cardinal of cofinality $\aleph_1$ such that $\l^{\aleph_1} > 
\l^{+(2^{\aleph_0})^+}$ 
then we have a counterexample as 
$(\l^{+\o})^{\aleph_0} < \l^{+(2^{\aleph_0})^+}$ (by \cite{ShA2, Ch\.~XIII, 
5\.1}).
 
The rest of this paper is devoted to the proof that the condition is 
necessary.
 
Assume that Tarski's conjecture fails, and let $\b$ be a limit ordinal for 
which there exists a sequence $\{\s_{\xi}\}_{\xi < \b}$ that gives a 
counterexample:
$$
\prod_{\xi < \b} \aleph_{\s_{\xi}} < \aleph_{\a}{}^{\kappa}, \tag 2
$$
where
$$
\kappa = |\b|\ \text{and}\ \a = \lim_{\xi < \b} \s_{\xi}.
$$
 
\proclaim{Lemma 1} If $(2)$ holds then $cf \b < \kappa < \b$, and there 
exists an ordinal $\g < \a$ such that $\aleph_{\g}{}^{\kappa} > \aleph_{\a}$.
\endproclaim
 
\demo{Proof} If (2) holds then $\b$ does not have $|\b|$ disjoint cofinal
subsets, and it follows that $\b$ is not a cardinal, and that cf $\b <
|\b|.$
 
Assuming that
$\aleph_{\g}{}^{\kappa} \le \aleph_{\a}$ 
holds for all $\g < 
\a$,  we pick a cofinal sequence $\{\a_i\}_{i < cf \b}$ with 
limit $\a$, and then
$$
\aleph_{\a}{}^{\kappa} =( \sum_{i < cf \b} 
\aleph_{\a_i})^{\kappa} \le \prod_{i < cf \b} 
\aleph_{\a_i}{}^{\kappa} \le \prod_{i < cf \b} \aleph_{\a} = 
\aleph_{\a}{}^{cf \b}= \prod_{i < cf \b} {\aleph_{\a_i}} \le
\prod_{\xi <\b} \aleph_{\s_\xi},
$$
contrary to $(2)$.
\qed
\enddemo
 
Now consider the shortest counterexample to Tarski's conjecture.
 
\proclaim{Lemma 2} If $\b$ is the least ordinal for which $(2)$ holds then $\b 
= \kappa + \o$ where $\kappa$ is an uncountable cardinal.
\endproclaim
 
\demo{Proof} Without loss of generality, the sequence $\s$ is continuous.
(We can replace each $\s_\xi$ by the limit of the sequence at $\xi$, for
each limit ordinal $\xi.$)
 
Let $\kappa = |\b|.$  We 
claim that for every limit ordinal $\eta 
< \b$,  $\aleph_{\s_{\eta}}{}^{\kappa} < 
\aleph_{\a}$.  If this were not true then, because $\b > \kappa$, there would 
be a limit ordinal $\eta$ 
such that $\kappa \le \eta
< \b$ and that 
$\aleph_{\s_{\eta}}{}^{|\s_{\eta}|} \ge \aleph_{\a}{}^{\kappa} 
> \prod_{\xi < \eta} 
\aleph_{\s_{\xi}}$, which would make the sequence $\{\s_{\xi}\}_{\xi < \eta}$ 
a counterexample to Tarski's conjecture as well, contrary to the minimality of 
$\b$.
 
Thus $\b = \d + \o$ for some limit ordinal $\d$.  It is clear that the 
sequence
$$
\{\aleph_{\s_{\xi}}: \xi \le \kappa\ \text{or}\ \xi > \d\}
$$
of length $\kappa + \o$ is also a counterexample, and by the minimality of 
$\b$ we have $\b = \kappa + \o$. \qed
\enddemo
 
Now consider the least ordinal $\g$ such that $\aleph_{\g}{}^{\kappa} > 
\aleph_{\a}$.  We shall show that $cf \g = \kappa$ (and so $\kappa$ is a 
regular uncountable cardinal).  We also establish other properties of 
$\aleph_{\g}$.
 
\proclaim{Lemma 3} If Tarski's conjecture fails, then there is a cardinal 
$\aleph_{\g}$ of uncountable cofinality $\kappa$ such that $\g > \kappa$, and 
that
$$
\align
{}&\text{for every $\nu < \g$, $\aleph_{\nu}{}^{\kappa} < \aleph_{\g}$} \tag 3 
\\ \vspace{1\jot}
{}&\aleph_{\g}{}^{\kappa} > \aleph_{\g+\o}{}^{\aleph_0}. \tag 4
\endalign
$$
\endproclaim
 
\demo{Proof} Let $\b = \kappa + \o$ be the least ordinal for which $(2)$ 
holds, for some increasing continuous sequence $\{\s_{\xi}: \xi < \b\}$ with 
limit $\a$, and let $\g$ be the least ordinal such that $\aleph_{\g}{}^{\kappa} 
> \aleph_{\a}$.
 
First we observe that for every $\nu < \g$, $\aleph_{\nu}{}^{\kappa} < 
\aleph_{\g}$.  This is because if $\aleph_{\nu}{}^{\kappa} \ge \aleph_{\g}$ 
then 
$\aleph_{\nu}{}^{\kappa} \ge \aleph_{\g}{}^{\kappa} > \aleph_{\a}$,  
contradicting the minimality of $\g$.
 
As a consequence, we have $cf \g \le \kappa$:  otherwise, we would have 
$\aleph_{\g}{}^{\kappa} = \sum_{\nu < \g} \aleph_{\nu}{}^{\kappa} = \aleph_{\g} < 
\aleph_{\a}$, a contradiction.  Also, if $\g = \lim_{i \to cf \g} 
\g_i$, then
$
\aleph_{\g}{}^{\kappa} = \left( \sum_{i < cf \g} 
\aleph_{\g_i}\right)^{\kappa} \le \prod_{i < cf \g} 
\aleph_{\g_i}{}^{\kappa} \le \prod_{i < cf \g} \aleph_{\g} = 
\aleph_{\g}{}^{cf \g}
$
and so we have $$\aleph_{\g}{}^{cf \g} = \aleph_{\g}{}^{\kappa}.$$
 
Since $\aleph_{\a} < \aleph_{\g}{}^{\kappa}$, we have $\aleph_{\a}{}^{\kappa} \le 
\aleph_{\g}{}^{\kappa} = \aleph_{\g}{}^{cf \g} \le 
\aleph_{\a}{}^{cf \g}$, and so $\aleph_{\a}{}^{cf \g} = 
\aleph_{\a}{}^{\kappa},$
and $\aleph_{\a}{}^{cf \g} > \prod_{\xi < \b} \aleph_{\s_{\xi}}$.  Hence the 
sequence
$$
\{\aleph_{\s_{\xi}}: \xi \le cf \g\ \text{or}\ \xi > 
\kappa\}
$$
of length $cf \g + \o$ is also a counterexample, and it follows that $\kappa 
= cf \g$.
 
For every limit $\eta < \b$ we have $\aleph_{\s_{\eta}}{}^{\kappa} < 
\aleph_{\a},$
and in particular
$\aleph_{\s_{\kappa}}{}^{\kappa} < \aleph_{\a}.$
Since $\aleph_{\g}{}^{\kappa} 
> \aleph_{\a}$, we have $\g > \kappa$.  Finally,
$$
\prod_{\xi < \b} \aleph_{\s_{\xi}} = \prod_{\xi < \kappa} \aleph_{\s_{\xi}} 
\cdot \prod_{n < \o} \aleph_{\s_{\kappa+n}} = \aleph_{\s_{\kappa}}{}^{\kappa} 
\cdot \aleph_{\a}{}^{\aleph_0} = \aleph_{\a}{}^{\aleph_0},
$$
and because $\aleph_{\g}{}^{\kappa} = \aleph_{\a}{}^{\kappa} > \prod_{\xi < \b} 
\aleph_{\s_{\xi}}$, we have $\aleph_{\g}{}^{\kappa} > \aleph_{\a}{}^{\aleph_0}$.  
Since $\a = {\underset {n \to \o} \to {\op{lim} \s_{\kappa+n}}} \ge \g + \o$, 
we have
$$
\aleph_{\g}{}^{\kappa} > \aleph_{\g+\o}{}^{\aleph_0},
$$
completing the proof. \qed
\enddemo
 
The cardinal $\aleph_{\g}$ obtained in Lemma~3 satisfies all the conditions 
stated in the Theorem except for the requirement that its cofinality be 
$\aleph_1$.  Thus the following lemma will complete the proof:
 
\proclaim{Lemma 4} Let $\aleph_{\g}$ be a singular cardinal of cofinality 
$\kappa > \aleph_1$ such that $\g > \kappa$ and that
$$
\text{for every $\nu < \g$, $\aleph_{\nu}{}^{\kappa} < \aleph_{\g}$.} \tag 5
$$
Assume further that for every $\d$, $\o_1 < \d < \g$, of cofinality 
$\aleph_1$,
$$
\text{if for every $\nu < \d$, $\aleph_{\nu}{}^{\aleph_1} < \aleph_{\d}$, then } 
\aleph_{\d}{}^{\aleph_1} \le \aleph_{\d+\o}{}^{\aleph_0} . \tag 6
$$
Then $\aleph_{\g}{}^{\kappa} \le \aleph_{\g+\o}{}^{\aleph_0}$.
\endproclaim
 
Lemma~4 implies that the least $\g$ in Lemma~3 has cofinality $\aleph_1$, and 
the theorem follows.  The rest of the paper is devoted to the proof of 
Lemma~4.  We use the second author's analysis of $pcf$.
 
\demo{Definition} If $A$ is a set of regular cardinals, let
$$
\prd A = \{f : \text{ dom} f = A \text{ and } f(\l) < \l \text{ for all }
\l \in A \}.
$$
If $I$ is an ideal on $A$ then $\prd A/I$ is a partially ordered set under
$$
f \le_I g \text{  iff  } \{\l : f(\l) > g(\l)\} \in I,
$$
and similarly for filters on $A.$ If $D$ is an ultrafilter on $A$, then
$\prd A/D$ is a linearly ordered set, and $cf(\prd A/D)$ denotes its
cofinality. Let
$$
pcf(A) = \{cf(\prd A/D): D\ \text{an ultrafilter on}\ A\}.
$$
It is clear that
$$
\align
A \subseteq pcf(A), A_1 \subseteq A_2\ \text{implies}\ &pcf(A_1) 
\subseteq pcf(A_2),\ \text{and} \\
&pcf(A_1 \cup A_2) = pcf(A_1) \cup pcf(A_2),
\endalign
$$
and it 
is not difficult to show (using ultrapowers of ultrapowers) that
$$
\align
\text{if}\ |pcf(A)| < \min A\ \text{then}\ &pcf(pcf(A)) = pcf(A)\ \text{and} 
\\
&pcf(A)\ \text{has a greatest element.}
\endalign
$$
\enddemo
 
\proclaim{Theorem {\rm (Shelah \cite{Sh345})}} If $2^{|A|} < \min(A)$ then there 
exists a family $\{B_{\nu}: \nu \in pcf(A)\}$ of subsets of $A$ such that
$$
\text{for every ultrafilter $D$ on $A$, $cf(\prd A/D) =$ the least $\nu$ such 
that $B_{\nu} \in D$.} \tag 7
$$
For every $\l \in pcf(A)$ there exists a family $\{f_{\a}: \a < \l\} \subseteq 
\prd A$ such that
$$
\aligned
{}&\text{$\a < \b$ implies $f_{\a} < f_{\b} \mod J_{<\l}$, where $J_{<\l}$ is the 
ideal generated} \\
{}&\text{by $\{B_{\nu}: \nu < \lambda\}$, and the $f_{\a}$'s are 
cofinal in $\prd B_{\lambda} \mod J_{< \lambda}$.} 
\endaligned \tag 8
$$ \qed
\endproclaim
 
An immediate consequence of $(7)$ is that $|pcf(A)| \le 2^{|A|}$.  The sets 
$B_{\nu}$  ($\nu \in pcf(A))$ are called {\it generators} for $A$.  Note that 
$\max B_{\nu} = \nu$ when $\nu \in A$, and that $\max (pcf(B_{\nu}))
= \nu$ for all $\nu$.
 
We shall use some properties of generators. 
 
\proclaim{Lemma 5 {\rm \cite{Sh345}}} 
Let $B_{\nu}$ be generators for $A$.
For every $X \subseteq A$ there exists a 
finite set $F \subseteq pcf(X)$ such that $X \subseteq \bigcup\{B_{\nu}: \nu \in 
F\}$.
\endproclaim
 
\demo{Proof} Let $Y = pcf(X)$, and assume that the lemma fails.  Then $\{X 
- B_{\nu}: \nu \in Y\}$ has the finite intersection property and so there is 
an ultrafilter $D$ on $A$ such that $X \in D$ and $B_{\nu} \notin D$ for all 
$\nu \in Y$.  Let $\mu = cf(\prd A/D)$.  Then $\mu \in pcf(X)$ and by $(7)$, 
$B_{\mu} \in D$.  A contradiction. \qed
\enddemo
 
For each $X \subseteq A$, let $s(X)$ (a {\it support} of $X$) denote a finite 
set $F \subseteq pcf(X)$ with the property that $X \subseteq \bigcup_{\nu 
\in F} B_{\nu}$.
 
The set $pcf(A)$ has a set of 
generators that satisfy a transitivity condition:
 
\proclaim{Lemma 6 {\rm \cite{Sh345}}} Assume that $2^{|A|} < \min(A)$ and let 
${\bar A} = pcf(A)$.  Then 
$pcf({\bar A})={\bar A}$ and
${\bar A}$ has a set of generators 
$\{B_{\nu}: \nu \in {\bar A}\}$ that satisfy, in addition to $(7)$,
$$
\text{if $\xi \in B_{\nu}$ then $B_{\xi} \subseteq B_{\nu}$.} \tag 9
$$ \qed
\endproclaim
 
We use the transitivity to prove the next lemma.
 
\proclaim{Lemma 7} Assume that $2^{|A|} < \min(A)$, let ${\bar A} = pcf(A)$, 
let $B_{\nu}$, $\nu \in {\bar A}$, be transitive generators for ${\bar A}$, 
and for each $X \subseteq {\bar A}$ let $s(X)$ be a support of $X$.  If $A = 
\bigcup_{i \in I} A_i$, then
$$
{\bar A} = \bigcup\left\{pcf(B_{\nu}): \nu \in pcf\left( \bigcup_{i \in I} 
s(pcf( A_i))\right)\right\}.
$$
\endproclaim
 
\proclaim{Corollary} $\max({\bar A}) = \max pcf \bigcup_{i \in I} s(pcf(A_i))$.
\endproclaim
 
\demo{[Proof of Corollary} Let $\l = \max({\bar A})$; $\l \in pcf(B_{\nu})$ for 
some $\nu$ in $pcf(\bigcup_i s(A_i))$.  Since $\max(pcf(B_{\nu})) 
= \nu$, we have 
$\l \le \nu$.]
\enddemo
 
\demo{Proof} Let $X = \bigcup_{i \in I} 
s(pcf(A_i))$ and $F = s(X)$.  We have $$A = 
\bigcup_{i \in I} A_i \subseteq \bigcup_{i \in I} pcf(A_i) \subseteq \bigcup_{i \in I} 
\bigcup\{B_{\xi}: \xi \in s(pcf(A_i))\} =$$
$$= \bigcup\{B_{\xi}: \xi \in X\} \subseteq 
\bigcup\{B_{\xi}: \xi \in \bigcup_{\nu \in F} B_{\nu}\} \subseteq \bigcup_{\nu \in F} 
B_{\nu}$$ (the last inclusion is a consequence of 
transitivity $(9)$).  Therefore $${\bar A} 
= pcf(A) \subseteq pcf(\bigcup_{\nu \in F} B_{\nu}) = \bigcup_{\nu \in F} 
pcf(B_{\nu}) \subseteq \bigcup\{pcf(B_{\nu}): \nu \in pcf(X)\}.$$ \qed
\enddemo
 
Toward the proof of Lemma~4, let $\{\g_i: i < \kappa\}$ be a continuous 
increasing sequence of limit ordinals of cofinality $< \kappa$, such that 
$\lim_{i \to \kappa} \g_i = \g$, $2^{\kappa} < \aleph_{\g_0}$, and that for 
all $i < \kappa$,
$$
\text{for all $\nu < \g_i$, $\aleph_{\nu}{}^{\kappa} < \aleph_{\g_i}$.} \tag 10
$$
 
\proclaim{Lemma 8} There is a closed unbounded set $C \subseteq \kappa$ such 
that for all $n = 1,2,\dots$,
$$
\max pcf(\{\aleph_{\g_i+n}: i \in C\}) \le \aleph_{\g+n}. \tag 11
$$
\endproclaim
 
\demo{Proof} We show that for each $n$ there exists a closed unbounded set 
$C_n \subseteq \kappa$ such that $\max pcf(\{\aleph_{\g_i+n}: i \in C_n\}) \le 
\aleph_{\g+n}$.  To prove this, let $n \ge 1$ be fixed and let $A = 
\{\aleph_{\g_i+n}: i < \kappa\}$.  Let $\l$ be the least element of $pcf(A)$ 
above $\aleph_{\g+n}$ (if there is none there is nothing to prove).  Let 
$\{B_{\nu}: \nu \in pcf(A)\}$ be subsets of $A$ that satisfy $(7)$, and let 
$\{S_{\nu}: \nu \in pcf(A)\}$ be the subsets of $\kappa$ such that $B_{\nu} = 
\{\aleph_{\g_i+n}: i \in S_{\nu}\}$.  It suffices to prove that the set 
$S_{\aleph_{\g+1}} \cup \dots \cup S_{\aleph_{\g+n}}$ contains a closed 
unbounded set.
 
Thus assume that the set $S = \kappa - (S_{\aleph_{\g+1}} \cup \dots \cup 
S_{\aleph_{\g+n}})$ is stationary.  Let $J_{<\l}$ be the ideal on $A$ 
generated by $\{B_{\nu}: \nu < \l\}$.  By Shelah's Theorem there exists a 
family $\{f_{\a}: \a < \l\}$ in $\prd A$ such that $\a < \b$ implies $f_{\a} 
< f_{\b} \mod J_{<\l}$.  Since all the sets $B_{\nu}$, $\nu < \aleph_{\g}$, 
are bounded, we get a family $\{g_{\a}: \a < \l\}$ of functions on $S$ such 
that $g_{\a}(i) < \aleph_{\g_i+n}$ for all $i \in S$, and such that $\a < \b$ 
implies that $g_{\a}(i) < g_{\b}(i)$ for eventually all $i \in S$.  This 
contradicts the results in \cite{GH} by which,
under the assumption (5), any family of almost disjoint 
functions in $\prod_{i \in S} \aleph_{\g_i+n}$ has size at most 
$\aleph_{\g+n}$. \qed
\enddemo
 
\demo{Proof of Lemma 4} Let $\g$ be a singular cardinal of cofinality $\kappa 
> \aleph_1$ that satisfies $(5)$ and $(6)$.  Let $\l$ be a regular cardinal 
such that $\aleph_{\g} < \l \le \aleph_{\g}{}^{\kappa}$.  We shall prove that 
$\l \le \aleph_{\g+\o}{}^{\aleph_0}$.
 
Let $\{\g_i: i < \kappa\}$ be an increasing continuous sequence that satisfies 
$(10)$, and let $C$ be a closed unbounded subset of $\kappa$ given by Lemma~8.  
Let
$$
S = \{i \in C: cf \g_i = \aleph_1\}.
$$
As $\kappa \ge \aleph_2$, $S$ is a stationary subset of $\kappa$.
\enddemo
 
\proclaim{Lemma 9} There exist regular cardinals $\l_i$, $i \in S$, such that 
for each $i \in S$, $\aleph_{\g_i} < \l_i \le \aleph_{\g_i}{}^{\aleph_1}$, and 
an ultrafilter $D$ on $S$ such that $cf(\prod_{i \in S} \l_i/D) = \l$.
\endproclaim
 
\demo{Proof} Let $I_0$ be the nonstationary ideal on $S$.  
There are $\l$ cofinal subsets $X$ of $\o_{\g}$ of size $|X|=\kappa.$
For every such set $X$, let $F_X \in \prod_{i \in C} [\aleph_{\g_i}]^{\le
\kappa}$ be the function defined by $F(i)=X\cap \o_{\g_i}.$ Then when
$X \ne Y$, $F_X$ and $F_Y$ are eventually distinct.
 
For every $i \in S$ we have $\aleph_{\g_i}{}^{\kappa}= \aleph_{\g_i}{}^
{\aleph_1}$ (by (10)), and so 
there exist $\l$ $I_0$-distinct functions in 
$\prod_{i \in S} \aleph_{\g_i}{}^{\aleph_1}$.  [$f$~and $g$ are 
$I_0$-{\it distinct} if $\{i: f(i) = g(i)\} \in I_0$.] 
 
Consider the partial 
ordering $f <_{I_0} g$ defined by $\{i: f(i) \ge g(i)\} \in I_0$; since $I_0$ 
is $\s$-complete, $<_{I_0}$ is well-founded.  Let $g$  
be a $<_{I_0}$-minimal function with 
the property that 
$g(i) \le \aleph_{\g_i}{}^{\aleph_1}$ and that there are
there are $\l$ $I_0$-distinct functions below $g.$
 
Let $I$ be the extension of $I_0$ generated by all the stationary subsets $X$ 
of $S$ that have the property that $g$ is not minimal on $I_0[X]$ (i\.e\. 
there is a function $g'$ such that $g'(i) < g(i)$ almost everywhere on $X$ and 
below $g'$ there are $\l$ $I_0$-distinct functions).
\enddemo
 
\demo{Claim} $I$ is a normal $\kappa$-complete ideal on $S$.
\enddemo
 
\demo{[Proof} Let $X_i$, $i < \kappa$, be sets in $I$, and let for each $i < 
\kappa$, $g_i < g$ on $X_i$ and $\<h_{\xi}^i: \xi < \l\>$ witness that $X_i 
\in I$.  Then one constructs witnesses ${\bar g}$ and $\<{\bar h}_{\xi}: \xi < 
\l\>$ for $X = \{j \in \kappa: j \in \bigcup_{i < j} X_i\}$ by letting ${\bar 
g}(j) = g_i(j)$ and ${\bar h}_{\xi}(j) = h_{\xi}^i(j)$ where $i$ is some $i < 
j$ such that $j \in X_i$.
 
For example, let us show that ${\bar h}_{\xi}$ and ${\bar h}_{\eta}$ are 
$I_0$-distinct if $\xi \ne \eta.$ Assume that ${\bar h}_{\xi}=
{\bar h}_{\eta}$ on a stationary subset $S_1$ of $S$. Then on a stationary
subset $S_2$ of $S_1$ the $i$ less than $j \in S_2$ chosen such that
$j \in X_i$ is the same $i,$ and we have $h_{\xi}^i =h_{\eta}^i$ on $S_2,$
a contradiction.]
\enddemo
 
Let $\{h_{\xi}: \xi < \l\}$ be a family of $I_0$-distinct functions below $g$.
 
\demo{Claim} For every $h <_I g$ there is some $\xi_0 < \l$ such that for all 
$\xi \ge \xi_0$, $h <_I h_{\xi}$.
\enddemo
 
\demo{[Proof} If there are $\l$ many $\xi$'s such that $h \ge h_{\xi}$ on an 
$I$-positive set, then (because $2^{\kappa} < \l$) there is an $I$-positive 
set $X$ such that $h \ge h_{\xi}$ on $X$ for $\l$ many $\xi$, but this 
contradicts the definition of $I$.]
\enddemo
 
Using this Claim, one can construct a $<_I$-increasing $\l$-sequence (a 
subsequence of $\{h_{\xi}: \xi < \l\}$) of functions that is $<_I$-cofinal in 
$\prod_{i \in S} g(i)$.  Let $\l_i = cf g(i)$, for each $i \in S$.  The 
product $\prod_{i \in S} \l_i$ has a $<_I$-cofinal $<_I$-increasing sequence of 
length $\l$, and since $I$ is a normal ideal, we have $\l_i > \aleph_
{\g_i}$ for 
$I$-almost all $i$.  Now if $D$ is any ultrafilter extending the dual of $I$, 
$D$ satisfies $cf(\prod_{i \in S} \l_i/D) = \l$. \qed
 
\medpagebreak
Back to the proof of Lemma~4.  For each $i \in S$ we have a regular cardinal 
$\l_i$ such that $\aleph_{\g_i} < \l_i \le \aleph_{\g_i}{}^{\aleph_1}$.  By the 
assumption $(6)$ we have $\aleph_{\g_i}{}^{\aleph_1} \le 
\aleph_{\g_i+\o}{}^{\aleph_0}$, and so $\l_i \le \aleph_{\g_i+\o}{}^{\aleph_0}$.  
We use the following result:
 
\proclaim{Theorem {\rm (Shelah \cite{ShA2}, Chapter XIII, 5\.1)}} Let 
$\aleph_{\d}$ be such that $\aleph_{\d}{}^{\aleph_0} < \aleph_{\d+\o}$.  Then 
for every regular cardinal $\mu$ such that $\aleph_{\d} < \mu \le 
\aleph_{\d+\o}{}^{\aleph_0}$ there is an ultrafilter $U$ on $\o$ such that 
$cf(\prod_{n \in \o} \aleph_{\d+n}/ U) = \mu$. \qed
\endproclaim
 
We apply the theorem to each $\aleph_{\g_i}$, and obtain for each $i \in S$ an 
ultrafilter $U_i$ on $\o$ such that $cf(\prod_{n \in \o} \aleph_{\g_i+n}/U_i) 
= \l_i$.  Combining the ultrafilters $U_i$ with the ultrafilter $D$ on 
$S$ from Lemma~9 we get an ultrafilter $U$ on the set
$$
A = \{\aleph_{\g_i+n}: i \in S,\,n = 1,2,\dots\}
$$
such that $cf(\prd A/U) = \l$.  Hence $\l \in pcf(A)$.
 
We shall now complete the proof of Lemma~4 by showing that $\max pcf(A) \le 
\aleph_{\g+\o}{}^{\aleph_0}$.
 
We have $A = \bigcup_{n=1}^{\i} A_n$, where
$$
A_n = \{\aleph_{\g_i+n}: i \in S\},
$$
and since $2^{|A|} = 2^{\kappa} < \min(A)$, we apply the corollary of
Lemma~7 and get
$$
\max pcf(A) = \max pcf \bigcup_{n=1}^{\i} s(pcf(A_n)),
$$
where for each $n$, $s(pcf(A_n))$ is a finite subset of $pcf(pcf(A_n))=
pcf(A_n)$.
 
Let $E = \bigcup_{n=1}^{\i} s(pcf(A_n))$.  Since (by Lemma~8) $\max pcf(A_n) \le 
\aleph_{\g+n}$ for each $n$, $E$ is a countable subset of $\aleph_{\g+\o}$.  
Hence $\max pcf(E) \le \aleph_{\g+\o}{}^{\aleph_0}$, and so
$$
\lambda \le \max pcf(A) = \max pcf(E) \le \aleph_{\g+\o}{}^{\aleph_0}.
$$
\qed
 
\bigskip
\heading References
\endheading
 
\bigskip
 
\ref \key[{\bf GH}]
\by F\. Galvin and A\. Hajnal
\paper Inequalities for cardinal powers
\jour Annals of Math\. \vol 101 \yr 1975 \pages 491--498
\endref
 
\ref \key[{\bf M}]
\by M\. Magidor
\paper On the singular cardinals problem {\rm I}
\jour Israel J\. Math\. \vol 28 \yr 1977 \pages 1--31
\endref
 
\ref \key[{\bf ShA2}]
\manyby S\. Shelah
\book Proper Forcing
\publ Springer--Verlag Lecture Notes 940 \yr 1982
\endref
 
\ref \key[{\bf Sh345}]
\bysame
\paper Products of regular cardinals and cardinals invariants of products of 
Boolean algebras
\jour Israel J\. Math\.
\vol 70 \yr 1990 \pages 129--187
\endref
 
\ref \key[{\bf Sh355}]
\bysame
\paper $\aleph_{\o+1}$ has a Jonsson algebra
\jour preprint
\endref
 
\ref \key[{\bf T}]
\by A\. Tarski
\paper Quelques th\'eor\`emes sur les alephs
\jour Fundamenta Mathematicae \vol 7 \yr 1925 \pages 1--14
\endref
 
\enddocument